\documentclass[12pt]{amsart}

\usepackage{amsmath}
\usepackage{amsthm}
\usepackage{amsfonts}
\usepackage{amssymb,amsmath,latexsym}
\usepackage{graphicx}

\theoremstyle{plain}
\newtheorem{theorem}{Theorem}[section]
\newtheorem{lemma}{Lemma}[section]

\theoremstyle{definition}
\newtheorem{remark}{Remark}[section]

\title{On Chordal and Bilateral SLE in multiply connected domains}

\author{Robert~O. Bauer}
\address{Department of Mathematics\\ 
	University of Illinois at Urbana-Champaign\\ 
	1409 West Green Street \\ 
	Urbana, IL 61801, USA}
\email{rbauer@math.uiuc.edu}
\thanks{The research of the first author was supported by 	NSA grant H98230-04-1-0039}
\author{Roland~M. Friedrich}
\address{Institute for Advanced Study\\ 
	Princeton, NJ 08540, USA}
\email{rolandf@ias.edu}
\thanks{The research of the second author was supported by 	NSF grant DMS-0111298.}

\begin{document}

\begin{abstract}
We discuss the possible candidates for conformally invariant random non-self-crossing curves which begin and end on the boundary of a multiply connected planar domain, and which satisfy a Markovian-type property. We consider both, the case when the curve connects a boundary component to itself (chordal), and the case when the curve connects two different boundary components (bilateral).  We establish  appropriate extensions of Loewner's equation to multiply connected domains for the two cases. We show that a curve in the domain induces a motion on the boundary and that this motion is enough to first recover the motion of the moduli of the domain and then, second, the curve in the interior. For random curves in the interior we show that the induced random motion on the boundary is not Markov if the domain is multiply connected, but that the random motion on the boundary together with the random motion of the moduli forms a Markov process. In the chordal case, we show that this Markov process satisfies Brownian scaling and discuss how this limits the possible conformally invariant random non-self-crossing curves. We show that the possible candidates are labeled by a real constant and a function homogeneous of degree minus one which describes the  interaction of the random curve with the boundary. We show that the random curve has the locality property if the interaction term vanishes and the real parameter equals six.
\end{abstract}

\maketitle

\section{Introduction} 

In this paper we discuss the possible candidates for a mathematically rigorous notion of conformally invariant random non-self-crossing curves which begin and end on the boundary of a multiply connected planar domain, and which satisfy a Markovian-type property. The Markovian type property means that the random curves can be developed dynamically as a (locally) growing family of random compacts. We aim to proceed in the spirit of Schramm, who deduced that, under an additional reflection symmetry, there is only a one parameter family of such random curves in simply connected domains, which he termed Stochastic Loewner Evolutions, see \cite{schramm:2000}. As such conformally invariant random growing compacts are conjectured to arise as scaling limits of interfaces of 2-dimensional statistical mechanical systems at criticality, Schramm had with one stroke identified what those limits can be. This has many consequences and applications, see \cite{lsw:2001a}, \cite{lsw:2003}, \cite{lsw:2004}, \cite{werner:2003}, and references therein.

Statistical mechanical systems have been studied in discrete approximations of multiply  connected domains and Riemann surfaces, see \cite{LPS}, and \cite{Ai}, and the connections with conformal field theory (CFT) indicate that the stochastic Loewner evolution should also extend to multiply connected domains and Riemann surfaces. 

For multiply connected domains the situation is already more subtle when compared to the simply connected case, because moduli spaces enter the picture and, as we will show, one has to consider interactions with these moduli. 

Families of random compacts from the boundary to the boundary now come in two flavors, as the random compact may grow to either connect a boundary component to itself (the chordal case) or it may grow to connect two different boundary components. We call the latter the {\it bilateral } case.  

The radial case, treated in \cite{BF:2004a}, where the random compact grows from the boundary to an interior point, can be considered  as  a limit of the bilateral case, when the boundary component the random compact grows towards shrinks to a point. This can be made precise, see \cite{komatu:1943}. 

Our procedure rests on an appropriate extension of Loewner's equation to the multiply connected case. In the simply connected case, Loewner's equation allows to encode a simple curve in a domain $D$ which has one endpoint on the boundary $\partial D$ by a continuous motion on the boundary, see \cite{loewner:1923}. In the multiply connected case, we show in Theorem \ref{T:CL} and Theorem \ref{T:BKL} that a simple curve induces a motion on the boundary of the domain. To recover the curve inside the domain requires also the knowledge of the moduli ${\bf M}$ (which describe the conformal equivalence class), as the curve grows. We show in Theorem \ref{T:Cmm} and Theorem \ref{T:Bmm} that these moduli can be recovered from the boundary motion and thus, once these moduli have been obtained, the curve in the interior itself.

A growing random non-self-crossing curve in a multiply connected domain can then also be encoded into a random motion $\xi(t)$ on the boundary. However, if the connectivity is greater than one, then $\xi$ cannot be Markov. We show in Section \ref{S:CSLE} that in the chordal case the boundary motion $\xi$ together with the motion of moduli ${\bf M}$ is a Markov process, and that it satisfies Brownian scaling. 

These facts dramatically reduce the number of possible diffusions. Indeed, in addition to a real parameter $\kappa$ one is only free to choose a function $A$ which is homogeneous of degree minus one in the variables $\xi$ and ${\bf M}$. The term $A$ measures the interaction of the random growing compact with the boundary (for example if it is desired that the random set avoids the interior boundary components). 

$\text{SLE}(\kappa,\rho)$, see \cite{lsw:2003}, and \cite{dubedat:2004}, also fits naturally into this framework. There, the random compact grows into the upper half-plane, the boundary is the real axis, and the interaction is with a finite number of points on the real axis and given in terms of the simplest  homogeneous function of degree minus one, $1/x$. Even though the upper half-plane is simply connected, the marked points on the boundary can serve as moduli and then $\text{SLE}(\kappa,\rho)$ is given by a particular moduli diffusion. 

For multiply connected domains it is natural to look for an interaction $A$ which is expressed in terms of domain functionals such as the Green function. Appropriate combinations of derivatives of the Green function are homogeneous of degree minus one in $\xi$ and the moduli. The `harmonic random Loewner chains' studied in \cite{zhan:2004} are a particular example of this. 

In our opinion the only further reduction in possible diffusions $(\xi,{\bf M})$ are regularity requirements on the homogenous function $A$. In particular we may wish to allow only functions which are analytic. We propose to call the growing family of random compacts obtained by solving the chordal Loewner equation \eqref{E:CL} for a diffusion $(\xi,{\bf M})$ associated to an analytic function $A$ homogenous of degree minus one by \eqref{E:syst},  {\em chordal stochastic Loewner evolution}.    

Finally, we would like to point out very briefly several physical aspects, which shall
be discussed in more detail in a separate publication. The classical Green function
is basically the two-point function of the bosonic free field with central charge $1$.
This elucidates the appearance of the number $-1$, the degree of homogeneity.
Further, scale invariance is related to the operator $L_0$, the infinitesimal generator of scale
transformations. It is an element, as is well known~\cite{BPZ}, of the Virasoro algebra, and in particular
an element of the sub-algebra corresponding to the group of global conformal transformations.

Interactions with boundary components such as those mentioned above, are conveniently modelled in the Coulomb gas formalism,
and correspond to insertions of different operators (i.e. currents or vertex operators) into the correlator.

As we are interested in describing the measure on random paths in a dynamical way, the moduli have to change, due to the deformation of the domain and/or the presence of marked points. This is in agreement with the global framework introduced in~\cite{FK, KBonn},  where correlators are
modelled as sections of a twisted version of a determinant line bundle over the appropriate moduli space.

In the case of the upper half plane, the CFT content of $\text{SLE}(\kappa,\rho)$ has been studied and explained in~\cite{cardy:2004}.

\section{Bilateral and Chordal standard domains}\label{S:canonical}




\subsection{Harmonic Measures}
Denote $D$ a region of connectivity $n>1$ in the complex plane. The components of the complement in the extended complex plane are denoted by $E_1, E_2,\dots, E_n$. We assume that no $E_k$ reduces to a point and that there is a unique unbounded component $E_n$. By applying preliminary conformal maps, we may assume that $D$ is bounded by an outer contour $C_n$ and $n-1$ inner contours $C_1,\dots, C_{n-1}$, where the contours are oriented such that $D$ lies to the left in the direction of the contour. Denote $\omega_k(z)$ the solution to the Dirichlet problem in $D$ with the boundary values $1$ on $C_k$ and $0$ on the other contours. We have $0<\omega_k(z)<1$ in $D$ and 
\begin{equation}\label{E:harmonic sum}
	\omega_1(z)+\omega_2(z)+\cdots+\omega_n(z)=1.
\end{equation}
$\omega_k(z)$ is called the harmonic measure of $C_k$ in $z$. The conjugate harmonic differential of $\omega_k$ has periods
\begin{equation}\label{E:omega-periods}
	\alpha_{kj}=\int_{C_j}*d\omega_k
	=\int_{C_j}\frac{\partial\omega_k}{\partial n}\ ds
\end{equation}
along $C_j$. Here, $\partial/\partial n$ denotes the normal derivative to the right of the direction of the contour, and $ds$ stands for arc-length measure. It is well known, \cite{nehari:1952}, that the $(n-1)\times(n-1)$ matrix $\boldsymbol{\alpha}$ with entries $\alpha_{kj}$, $1\le k,j\le n-1$, is positive definite and symmetric. In particular, the linear system
\begin{align}\label{E:one period}
	\lambda_1 \alpha_{11}+\lambda_2\alpha_{21}+\cdots
	+\lambda_{n-1}\alpha_{n-1, 1}&=2\pi\notag\\
	\lambda_1 \alpha_{12}+\lambda_2\alpha_{22}+\cdots
	+\lambda_{n-1}\alpha_{n-1, 2}&=0\notag\\
	&\ \vdots\\
	\lambda_1 \alpha_{1, n-1}+\lambda_2\alpha_{2,n-1}+\cdots
	+\lambda_{n-1}\alpha_{n-1, n-1}&=0\notag
\end{align}
has a unique solution. It follows from \eqref{E:harmonic sum} that any solution of \eqref{E:one period} also solves
\[
	\lambda_1\alpha_{1n}+\lambda_2\alpha_{2n}+\cdots
	+\lambda_{n-1}\alpha_{n-1,n}=-2\pi.
\]
Thus there is a multiple-valued integral $F(z)$ with periods $\pm2\pi i$ along $C_1$ and $C_n$ and all other periods equal to zero, the real part being constant equal to $\lambda_k$ on $C_k$ (we set $\lambda_n=0$). The function $f(z)=e^{-F(z)}$ is then single-valued and one can show, \cite{ahlfors:1966}, that $f$ maps $D$ conformally onto the annulus $e^{-\lambda_1}<|w|<1$ minus $n-2$ concentric arcs situated on the circles $|w|=e^{-\lambda_k}$, $k=2,\dots,n-1$. We call such a circularly slit annulus a {\em bilateral standard domain}. By adding an imaginary constant to $F(z)$ we obtain another map onto a bilateral standard domain and we may normalize the map $f$ by requiring $f(z_0)=e^{-\lambda_1}$ for some $z_0\in C_1$. With this normalization we call $f$ the {\em canonical map} for $(D, z_0, C_n)$.


\subsection{Green Function}\label{ss:Green}
Denote $D$ again a region of finite connectivity which is bounded by contours $C_1,\dots, C_n$; this time the case $n=1$ is included.

We consider a point  $z_0\in D$ and solve the Dirichlet problem in $D$ with the boundary values $\ln|\zeta-z_0|$. The solution is denoted by $h(z)$. The function 
\[
	G(z)=G_D(z,z_0)=h(z)-\ln|z-z_0|
\]
is the Green function in $D$ with pole at $z_0$. It is the unique function which is harmonic in $D$ except at $z_0$, where it differs from $\ln|z-z_0|$ by a harmonic function, and which vanishes on the boundary of $D$. The Green function is conformally invariant in the sense that if $f:D\to D'$ is conformal, then
\begin{equation}\label{E:conformal Green}
	G_D(z,z_0)=G_{D'}(f(z),f(z_0)).
\end{equation}
The conjugate harmonic function of $G(z,z_0)$ is multiple-valued. It has the period $2\pi$ along a small circle about $z_0$, and the periods
\[
	p_k(z_0)=\int_{C_k}*dG(z,z_0),\quad k=1,\dots,n.
\]
It can be shown that $p_k(z_0)=-2\pi\omega_k(z_0)$, \cite{ahlfors:1966}. Let now $z_0\in C_n$. By linearity, $u(z)=\partial G(z,z_0)/\partial n_{z_0}$ is a harmonic function in $z$. Its conjugate differential has periods
\begin{align}
	A_k(z_0)&=\int_{C_k}*d u=\frac{\partial}{\partial n_{z_0}} 
	\int_{C_k}*d_{z} G(z,z_0)\notag\\
	&=-2\pi\frac{\partial}{\partial n_{z_0}}\omega_k(z_0).
\end{align}
Thus the linear combination $u+\lambda_1 \omega_1+\cdots+\lambda_{n-1}\omega_{n-1}$ is free from periods provided that 
\begin{equation}\label{E:no periods}
	\lambda_1\alpha_{1k}+\lambda_2\alpha_{2k}+\cdots
	+\lambda_{n-1}\alpha_{n-1,k}=-A_k,\quad k=1,\dots,n-1.
\end{equation}	
If we write ${\bf P}$ for the matrix $\boldsymbol\alpha/2\pi$, $\boldsymbol\lambda^T=(\lambda_1,\dots,\lambda_{n-1})$, and 
\[
	\partial\boldsymbol \omega(z_0)^T/\partial n
	=(\partial\omega_1(z_0)/\partial n,\dots, 
	\partial\omega_{n-1}(z_0)/\partial n),
\]
then the solution to \eqref{E:no periods} is given by
\[
	\boldsymbol\lambda
	={\bf P}^{-1}\frac{\partial\boldsymbol\omega(z_0)}{\partial n}.
\]
Hence
\[
	-i\left(\frac{\partial G(z,z_0)}{\partial n_{z_0}}
	+\boldsymbol\omega(z)^T
	{\bf P}^{-1}\frac{\partial\boldsymbol\omega(z_0)}{\partial n}\right)
\]
is the imaginary part of a single-valued analytic function $\Psi(z)$. It can be shown that $\Psi$ maps $D$ conformally onto the upper half-plane $\Im(w)>0$ minus  $n-1$ horizontal slits with imaginary parts 
\[
	\Im(w) 
	=-[{\bf P}^{-1}\partial\boldsymbol\omega(z_0)/\partial n]_j, 
	\quad j=1,\dots, n-1.
\]
Under this map, $\Psi(C_n)=\mathbb R$, and $\Psi(z_0)=\infty$. We call the upper half-plane minus a finite number of horizontal slits a {\em chordal standard domain}. If $D$ is contained in the upper half-plane and for some $x>0$ we have $(\mathbb R\backslash [-x,x])\subset C_n$ and $\zeta=\infty$, then, by adding an appropriate real constant, we may assume that $g\equiv\Psi/2$ satisfies the hydrodynamic normalization at infinity,
\begin{equation}\label{E:hydro-norm}
	\lim_{z\to\infty}(g(z)-z)=0.
\end{equation}
With this normalization, we call $g$ the canonical mapping for $D$.


\section{Evolution of slit mappings}

\subsection{Chordal Loewner equation}
Consider a chordal standard domain $D$. Let $\gamma:[0,t_{\gamma}]\to\overline{D}$ be a Jordan arc such that $\gamma(0)\in\mathbb R$, and $\gamma(0,t_{\gamma}]\subset D$. Let $g_t$ be the canonical mapping from $D\backslash\gamma[0,t]$ with the normalization \eqref{E:hydro-norm}, and denote $D_t$ the chordal standard domain $g_t(D\backslash\gamma[0,t])$. It is well known, see \cite{nehari:1952}, that $g_t$ solves the extremal problem 
\[
	a_1=\max
\]
among all univalent functions on $D\backslash\gamma[0,t]$ with expansion
\[
	z+\frac{a_1}{z}+\frac{a_2}{z^2}+\cdots,\quad
	a_k\in\mathbb R,
\]
near infinity. In particular, if $g_t(z)=z+a_t/z+o(1/|z|)$, then $a_{t^*}\le a_t$ whenever $0<t^*<t<t_{\gamma}$. In fact, a simple argument shows that 
\begin{equation}\label{E:strictly monotone}
	a_{t^*}< a_t\quad\text{if}\quad t^*<t.
\end{equation}
Thus we may assume that $a_t=2t$. We wish to find a differential equation for the family $\{g_t:t\in[0,t_{\gamma}]\}$. 

Denote $C_j(t), j=1,\dots,n$, the boundary components of $D_t$. We always have $C_n(t)=\mathbb R$. For $j=1,\dots, n-1$, let $y_j(t)$ be the imaginary part of (points on) the slit $C_j(t)$. Denote $\xi(t)$ the starting point on $\mathbb R$ of the Jordan arc $g_t(\gamma[t,t_{\gamma}])$ in $D_t$, i.e. $g_t(\gamma_t)$. For $0<t^*<t<t_{\gamma}$, set 
\[
	g_{t,t^*}=g_{t^*}\circ g_t^{-1}.
\]  
Then $g_{t,t^*}$ is a conformal map from $D_t$ onto $D_{t^*}\backslash g_{t^*}(\gamma[t^*,t])$. The point $\xi(t^*)=g_{t^*}(\gamma_{t^*})$ corresponds to two prime ends in  $D_{t^*}\backslash g_{t^*}(\gamma[t^*,t])$. Denote $\beta_0(t,t^*)$ and $\beta_1(t,t^*)$, with $\beta_0(t,t^*)<\beta_1(t,t^*)$, the pre-images of these prime ends under $g_{t,t^*}$, i.e.
\[
	g_{t,t^*}(\beta_0(t,t^*))=g_{t,t^*}(\beta_1(t,t^*))
	=g_{t^*}(\gamma_{t^*}).
\]
Then, if $x\in\mathbb R\backslash[\beta_0(t,t^*),\beta_1(t,t^*)]$, 
\[
	g_{t,t^*}(x)\in\mathbb R.
\]

Consider the analytic function 
\[
	z\mapsto g_{t,t^*}(z)-z,
\]
which satisfies
\begin{equation}\label{E:difference}
	g_{t,t^*}(z)-z=\frac{2(t^*-t)}{z}+o(1/|z|),
\end{equation}
and note that $z\mapsto\Im(g_{t,t^*}(z)-z)$ is harmonic and constant on each boundary component. By Poisson's formula 
\begin{equation}\label{E:poisson1}
	\Im(g_{t,t^*}(z)-z)=-\frac{1}{2\pi}\int_{\partial D_t}
	\Im(g_{t,t^*}(\zeta)-\zeta)\frac{\partial G_t(\zeta,z)}{\partial n_1}\ ds,
\end{equation}
where $G_t(\zeta,z)$ is the Green function for $D_t$ with pole at $z$. Note that there is no problem with integrability in \eqref{E:poisson1} because
\[
	\Im(g_{t,t^*}(\zeta)-\zeta)=\frac{-y}{x^2+y^2}+O(1/|\zeta|^2),
	\quad \zeta=x+iy,
\]
and 
\begin{equation}\label{E:imbound}
	\sup\{\Im(\zeta):\zeta\in\partial D_t\}<\infty.
\end{equation}
	Since $\Im(g_{t,t^*}(z)-z)$ has a single-valued harmonic conjugate, it is orthogonal to the real part of any Abelian differential of the first kind, see \cite{bauer.friedrich:2004}, and we have 
\begin{align}\label{E:poisson2}
	\Im&(g_{t,t^*}(z)-z)\notag\\
	&=-\frac{1}{2\pi}\int_{\partial D_t}
	\Im(g_{t,t^*}(\zeta)-\zeta)
	\left(\frac{\partial G_t(\zeta,z)}{\partial n_1}
	+\boldsymbol{\omega}_t(z)^T\boldsymbol{P}_t^{-1}\frac{\partial\boldsymbol{\omega}_t(\zeta)}{\partial n}\right) ds.
\end{align}
It follows from Section \ref{S:canonical} that
\[
	\int_{C_k(t)} \left(\frac{\partial G_t(\zeta,z)}{\partial n_1}
	+\boldsymbol{\omega}_t(z)^T\boldsymbol{P}_t^{-1}\frac{\partial\boldsymbol{\omega}_t(\zeta)}{\partial n}\right) ds=0,
	\quad k=1,\dots, n-1,
\]
and also that 
\[
	z\mapsto -i\left(\frac{\partial G_t(\zeta,z)}{\partial n_1}
	+\boldsymbol{\omega}_t(z)^T\boldsymbol{P}_t^{-1}\frac{\partial\boldsymbol{\omega}_t(\zeta)}{\partial n}\right)
\]
is the imaginary part of a single-valued analytic function $\Psi_t(z)=\Psi_t(z,\zeta)$. Thus, since  $\Im(g_{t,t^*}(\zeta)-\zeta)$ is constant on each $C_k(t)$, $k=1,\dots, n-1$, and identically zero on $\mathbb R\backslash[\beta_0(t,t^*),\beta_1(t,t^*)]$,
\begin{equation}\label{E:complex poisson2}
	g_{t,t^*}(z)-z
	=\frac{1}{2\pi}\int_{\beta_0(t,t^*)}^{\beta_1(t,t^*)}
	\Im(g_{t,t^*}(\zeta)-\zeta)\Psi_t(z,\zeta)\ d\zeta+ic,
\end{equation}
where $c$ is a real constant. Note that if $z\mapsto\tilde{\Psi}_t(z,\zeta)$ is another analytic function with the same imaginary part as $\Psi_t$, then 
\[
	\Psi_t(z,\zeta)-\tilde{\Psi}_t(z,\zeta)= a(\zeta),
\]
where $a$ is real and depends only on $\zeta$. We fix a normalization by requiring that 
\begin{equation}\label{E:normalPsi}
	\lim_{z\to\infty} \Psi_t(z,\zeta)=0.
\end{equation}
If we let $z\to\infty$, then $g_{t,t^*}(z)-z\to0$. By bounded convergence, the integral in \eqref{E:complex poisson2} converges to zero as well and it follows that $c=0$.  
Next, 
\[
	2(t^*-t)=\lim_{z\to \infty} z(g_{t,t^*}(z)-z)=f(0),
\]
where 
\[
	w\mapsto f(w)\equiv\frac{1}{w}\left[g_{t,t^*}\left(1/w\right)
	-\left(1/w\right)\right]
\]
is regular near zero. By the Schwarz reflection principle $g_{t,t^*}$ extends to the entire complex plane minus the slits $C_1(t),\dots, C_{n-1}(t)$, their conjugates, and the real interval $[\beta_0(t,t^*),\beta_1(t,t^*)]$. Denote $C$ the collection of these $2n-1$ finite slits. Then $f$ also extends to a corresponding domain with boundary $\tilde{C}$. From Cauchy's integral formula we have
\begin{align}
	2(t^*-t)&=\frac{1}{2\pi i}\int_{\tilde{C}}\frac{f(\zeta)}{\zeta}\ d\zeta
	=\frac{1}{2\pi i}\int_{\tilde{C}}
	\frac{g_{t,t^*}(1/\zeta)-1/\zeta}{\zeta^2}\ d\zeta\notag\\
	&=-\frac{1}{2\pi i}\int_C(g_{t,t^*}(\eta)-\eta)\ d\eta=-\frac{1}{2\pi }\int_C\Im(g_{t,t^*}(\eta)-\eta)\ d\eta,
\end{align}
where the final equality uses the fact that $d\eta$ is real for horizontal slits.
The slits $C_1(t),\dots,C_{n-1}(t)$ and their conjugates do not contribute to the last integral since $\Im(g_{t,t^*}(\eta)-\eta)$ takes the same value on both ``sides" of a given slit. For  the slit $[\beta_0(t,t^*),\beta_1(t,t^*)]$, $\Im(g_{t,t^*}(\eta)-\eta)$ takes opposite values on the upper and lower ``side" of the slit and, since the direction of integration is reversed, we finally get
\begin{equation}\label{E:integralrep}
	t^*-t=-\frac{1}{2\pi}\int_{\beta_0(t,t^*)}^{\beta_1(t,t^*)}
	\Im(g_{t,t^*}(\eta))\ d\eta.
\end{equation}

Setting $z=g_t(w)$ in \eqref{E:complex poisson2} we have
\begin{equation}\label{E:difference2}
	g_{t^*}(w)-g_t(w)\notag\\
	=\frac{1}{2\pi}\int_{\beta_0(t,t^*)}^{\beta_1(t,t^*)}
	\Im(g_{t,t^*}(\eta))\Psi_t(z,\eta) d\eta.
\end{equation}	
We are now ready to let $t^*\nearrow t$ in \eqref{E:difference2}. Note first that, for $\eta\in[\beta_0(t,t^*),\beta_1(t,t^*)]$, $\eta\mapsto\Im(g_{t,t^*}(\eta))$ is continuous and non-negative and that also $\eta\mapsto A(\eta):=\Psi_t(z,\eta)$
is continuous. Thus it follows from the mean-value theorem of integration and \eqref{E:integralrep} that 
\begin{align}
	\frac{1}{2\pi}&\int_{\beta_0(t,t^*)}^{\beta_1(t,t^*)}
	\Im(g_{t,t^*}(\eta)) A(\eta)\ d\eta\notag\\
	&=\frac{\Re(A(\eta'))+i\Im(A(\eta''))}{2\pi}
	\int_{\beta_0(t,t^*)}^{\beta_1(t,t^*)}
	\Im(g_{t,t^*}(\eta)\ d\eta\notag\\
	&=-\left[\Re(A(\eta'))+i\Im(A(\eta''))\right](t^*-t),
\end{align}
for some $\eta',\eta''\in[\beta_0(t,t^*),\beta_1(t,t^*)]$. Hence
\begin{equation}
	\lim_{t^*\nearrow t}\frac{ g_{t^*}(w)- g_t(w)}{t^*-t}\notag\\
	=-\Psi_t(z,\xi(t)).
\end{equation}
By the same argument we may let $t\searrow t^*$. On the right-hand side above we then only need to change $t$ to $t^*$ and introduce an overall minus sign. Thus we have established the following

\begin{theorem}[Chordal Loewner equation]\label{T:CL}
If $\gamma$ is a Jordan arc in a standard domain $D$ starting on $\mathbb R$ with the parametrization from above, and if $g_t$ is the canonical map for $D\backslash\gamma[0,t]$, then, using the notation from above, the family $\{g_t:t\in[0,t_{\gamma}]\}$ satisfies the equation
\begin{equation}\label{E:CL}
	\partial_t g_t(z)
	=-\Psi_t(g_t(z),\xi(t)),
\end{equation} 
with initial condition $g_0(z)=z$, and where $\Psi_t(z,\zeta)$ is the analytic function in $z$ with imaginary part
\[
	-\frac{\partial G_t(z,\zeta)}{\partial n_{\zeta}}
	-\boldsymbol{\omega}_t(z)^T\boldsymbol{P}_t^{-1}
	\frac{\partial\boldsymbol{\omega}_t(\zeta)}{\partial n},
\]
normalized by $\lim_{z\to\infty}\Psi_t(z,\zeta)=0$.
\end{theorem}

\begin{remark}
In the simply connected case, when $D=\mathbb H$ is the upper half-plane, the Green function is given by
\[
	G(z,w)=\Re\left(\ln\frac{z-\overline{w}}{z-w}\right).
\]
Thus, if $w=x+iy$,
\[
	-\Psi(z,w)=-i\frac{\partial}{\partial y}|_{y=0}\ln\frac{z-x+iy}{z-x-iy}=\frac{2}{z-w},
\]
and \eqref{E:CL} reduces to the well known chordal Loewner equation.
\end{remark}


\subsection{Bilateral Komatu-Loewner equation} The evolution of slit mappings in multiply connected domains was first studied by Komatu in \cite{komatu:1943} for the doubly connected case, and in \cite{komatu:1950} for general finite connectivity. Komatu treated this case by considering circular slit annuli. 

Consider a bilateral standard domain $D$ with inner radius $Q$. Let $\gamma:[0,t_{\gamma}]\to\overline{D}$ be a Jordan arc such that $\gamma(0)\in S^1$, and $\gamma(0,t_{\gamma}]\subset D$. Let $f_t$ be the canonical mapping from $D\backslash\gamma[0,t]$ with the normalization $f_t(Q)>0$, and denote $D_t$ the chordal standard domain $f_t(D\backslash\gamma[0,t])$. If $Q_t=f_t(Q)$, then it can be shown that $t\in[0,t_{\gamma}]\mapsto Q_t\in[Q,1]$ is continuous and strictly increasing, \cite{komatu:1950}. Thus we may assume that $\gamma$ is parametrized such that $t=\ln Q_t$. For this parameter it is shown in \cite{komatu:1950} that $t\mapsto f_t(z)$ is differentiable. An expression for the derivative is also given. However, the expression given there is not explicit enough for the purposes we have in mind. In particular, we will need to know that the vector field is itself a Lipschitz function in the moduli of the domain. We sketch a proof of what we call the bilateral Komatu-Loewner equation, leading to an expression of the derivative $\partial_t f_t$ in terms of the Green function, harmonic measures, their derivatives and harmonic conjugates. The argument is similar to the radial case, \cite{BF:2004a}. In fact, the radial case can be obtained as a limiting case from the bilateral case when $Q\to0$, \cite{komatu:1943}.

Denote $C_j(t), j=1,\dots,n$, the boundary components of $D_t$. We always have $C_n(t)=S^1$, and $C_1(t)=\{|z|=e^t\}$. For $j=2,\dots, n-1$, let $m_j(t)$ be the radial distance of the circular slit $C_j(t)$ from the origin. Denote $\xi(t)$ the starting point on $S^1$ of the Jordan arc $g_t(\gamma[t,t_{\gamma}])$ in $D_t$, i.e. $g_t(\gamma_t)$. For $\ln Q<t^*<t<t_{\gamma}\le0$, set 
\[
	g_{t,t^*}=g_{t^*}\circ g_t^{-1}.
\]  
Then $g_{t,t^*}$ is a conformal map from $D_t$ onto $D_{t^*}\backslash g_{t^*}(\gamma[t^*,t])$. The point $\xi(t^*)=g_{t^*}(\gamma_{t^*})$ corresponds to two prime ends in  $D_{t^*}\backslash g_{t^*}(\gamma[t^*,t])$. Denote $\exp(i\beta_0(t,t^*))$ and $\exp(i\beta_1(t,t^*))$, with $\beta_0(t,t^*)<\beta_1(t,t^*)$, the pre-images of these prime ends under $g_{t,t^*}$, i.e.
\[
	g_{t,t^*}(\exp(i\beta_0(t,t^*)))=g_{t,t^*}(\exp(i\beta_1(t,t^*)))
	=g_{t^*}(\gamma_{t^*}).
\]
Then, if $|z|=1$ and $\beta_1(t,t^*)\le\arg z\le\beta_0(t,t^*)+2\pi$, 
\[
	|g_{t,t^*}(z)|=1.
\]

The function 
\[
	z\mapsto\ln\frac{g_{t,t^*}(z)}{z}
\]
is analytic and single-valued throughout $D_t$. By Poisson's formula
\begin{equation}\label{E:poisson}
	\ln\left|\frac{g_{t,t^*}(z)}{z}\right|=-\frac{1}{2\pi}\int_{\partial D_t}
	\ln\left|\frac{g_{t,t^*}(\zeta)}{\zeta}\right|\frac{\partial G_t(\zeta,z)}{\partial n_1}\ ds,
\end{equation}
where $G_t(\zeta,z)$ is the Green function for $D_t$ with pole at $z$. Using orthogonality and the period relations as we did in the chordal case, it follows that
\begin{equation}\label{E:complex poisson3}
	\ln\frac{g_{t,t^*}(z)}{z}
	=-\frac{i}{2\pi}\int_{\beta_0(t,t^*)}^{\beta_1(t,t^*)}
	\ln\left|\frac{g_{t,t^*}(\zeta)}{\zeta}\right|\Psi_t(z,\zeta)\ ds +ic,
\end{equation}
for some real constant $c$. To eliminate $c$, we evaluate the identity \eqref{E:complex poisson3} at $z=q=e^t$ and then take the difference:
\begin{equation}\label{E:rep}
 	\ln\frac{g_{t,t^*}(z)}{z}-\ln\frac{q^*}{q}
	=-\frac{i}{2\pi}\int_{\beta_0(t,t^*)}^{\beta_1(t,t^*)}
	\ln\left|\frac{g_{t,t^*}(\zeta)}{\zeta}\right|\left[
	\Psi_t(z,\zeta)-\Psi_t(q,\zeta)\right] ds.
\end{equation}
By Cauchy's integral formula,
\begin{equation}\label{E:cauchy}
	0=\frac{1}{2\pi i}\int_{\partial D_t}
	\ln\left(\frac{g_{t,t^*}(\zeta)}{\zeta}\right)\frac{d\zeta}{\zeta}.
\end{equation}
In particular, the right-hand side of \eqref{E:cauchy} is real.  Since all boundary components are concentric circular arcs, $d\zeta/\zeta$ is purely imaginary along $\partial D_t$, i.e.
\[
	\frac{d\zeta}{\zeta}=i\ d\arg\zeta,\quad\zeta\in\partial D_t.
\]
Hence
\begin{align}
	0&=\frac{1}{2\pi}\int_{\partial D_t}
	\ln\left|\frac{g_{t,t^*}(\zeta)}{\zeta}\right|\ d\arg\zeta\notag\\
	&=\frac{1}{2\pi}\int_{\beta_0(t,t^*)}^{\beta_1(t,t^*)}
	\ln\left|g_{t,t^*}(e^{i\varphi})\right|\ d\varphi
	-\frac{1}{2\pi}\int_0^{2\pi}\ln\frac{q^*}{q}\ d\varphi\notag\\ 
	&\quad+\frac{1}{2\pi}\sum_{j=2}^{n-1}
	\int_{C_j(t)}\ln\frac{m_j(t^*)}{m_j(t)}\ d\arg\zeta.
\end{align}
Since the two ``sides" of $C_j(t)$ make opposite contributions,
\[
	\int_{C_j(t)}d\arg\zeta=0,\quad j=2,\dots, n-1,
\]
and we finally get
\begin{equation}\label{E:increment}
	t^*-t=\frac{1}{2\pi}\int_{\beta_0(t,t^*)}^{\beta_1(t,t^*)}
	\ln\left|g_{t,t^*}(e^{i\varphi})\right|\ d\varphi.
\end{equation}

Letting $z=g_t(w)$ in \eqref{E:rep}, we have
\begin{align}\label{E:difference3}
	\ln\frac{g_{t^*}(w)}{g_t(w)}&-(t^*-t)\notag\\
	&=-\frac{i}{2\pi}\int_{\beta_0(t,t^*)}^{\beta_1(t,t^*)}
	\ln|g_{t,t^*}(e^{i\varphi})| \left[
	\Psi_t(z,e^{i\varphi})-\Psi_t(q,e^{i\varphi})\right]ds.
\end{align}	
 
We now wish to let $t^*\nearrow t$ in \eqref{E:difference3}. Note first that, for $\varphi\in[0,2\pi]$, $\varphi\mapsto\ln|g_{t,t^*}(e^{i\varphi})|$ is continuous and non-positive and that also
\[
	\varphi\mapsto A(\varphi)
	:=\Psi_t(z,e^{i\varphi})-\Psi_t(q,e^{i\varphi})
\]
is continuous. Thus it follows from the mean-value theorem of integration that 
\begin{align}
	\frac{1}{2\pi(t^*-t)}&\int_{\beta_0(t,t^*)}^{\beta_1(t,t^*)}
	\ln\left|g_{t,t^*}(e^{i\varphi})\right| A(\varphi)\ d\varphi\notag\\
	&=\frac{\Re(A(\varphi'))+i\Im(A(\varphi''))}{2\pi(t^*-t)}
	\int_{\beta_0(t,t^*)}^{\beta_1(t,t^*)}
	\ln\left|g_{t,t^*}(e^{i\varphi})\right|\ d\varphi\notag\\
	&=\Re(A(\varphi'))+i\Im(A(\varphi'')),
\end{align}
for some $\varphi',\varphi''\in[\beta_0(t,t^*),\beta_1(t,t^*)]$. Hence
\begin{equation}
	\lim_{t^*\nearrow t}\frac{\ln g_{t^*}(w)-\ln g_t(w)}{t^*-t}
	=1+i[\Psi_t(z,\xi_t)-\Psi_t(e^t,\xi(t))].
\end{equation}
By the same argument we may let $t\searrow t^*$. On the right-hand side above we then only need to change $t$ to $t^*$ and introduce an overall minus sign. Thus we have established the following

\begin{theorem}[Bilateral Komatu-Loewner equation]\label{T:BKL}
If $\gamma$ is a Jordan arc in a standard domain $D$ starting on $S^1$ with the parametrization from above, and if $g_t$ is the canonical map for $D\backslash\gamma[0,t]$, then, using the notation from above, the family $\{g_t:t\in[\ln Q,t_{\gamma}]\}$ satisfies the equation
\begin{equation}\label{E:BKL}
\partial_t\ln g_t(z)
	=1+i[\Psi_t(g_t(z),\xi_t)-\Psi_t(e^t,\xi(t))],
\end{equation} 
with initial condition $g_{\ln Q}(z)=z$.
\end{theorem}


\section{Motion of moduli}\label{S:mm}

\subsection{Chordal case}
The right-hand side of the chordal Loewner equation, at time $t$, involves the Green function of the domain $D_t$, and also various functions derived from the Green function. Consequently, it does not make sense to ask for the solution of \eqref{E:CL} for a given continuous curve $t\mapsto\xi(t)$, since the vector-field on the right-hand side of \eqref{E:CL} is not specified by giving that information alone. To specify the Green function of $D_t$ we also need the {\em moduli} of the domain $D_t$. We will now consider what the appropriate moduli space is for our purposes and find a system of equations these moduli satisfy. Once this system is found, we can solve it for a given input $t\mapsto\xi(t)$, and then, in a second step, solve the radial Komatu-Loewner equation using $\xi$ and the moduli.

The geometric description of $D_t$ requires $3n-3$ real parameters, three for each (interior) slit, given, for example, by the imaginary components of the slits, i.e $y_j(t)$, $j=1,\dots, n-1$, and the real components
\[
	x_j(t)<x_j'(t),\quad j=1,\dots,n-1,
\]
determining the endpoints of the slit $C_j(t)$, $j=1\dots,n-1$. On the other hand, it is well known that two $n$-connected domains with non-degenerate boundary continua are conformally equivalent if $3n-6$ real parameters agree for $n>2$. If $n=2$ then there is only one real parameter describing the conformal class, and if $n=1$, then all such domains are conformally equivalent. 

The slits we wish to grow mark two points on one of the boundary continua, the beginning ($t=0$) and end point $(t=\infty)$ of the slit. Any $n$-connected planar domain with two marked boundary points on one boundary component is conformally equivalent to the upper half-plane with $n-1$ horizontal slits and such that the marked boundary points are mapped to $0$ and $\infty$. However, there is a one-parameter group of automorphisms, namely multiplication by $a>0$, which maps the slit upper half-plane onto a slit upper half-plane, while fixing 0 and $\infty$. It is now easy to see that the moduli space of $n$-connected planar domains with two marked boundary points on one of the boundary components  is $3n-4$ dimensional for all $n\ge2$, and zero dimensional if $n=1$. Nonetheless, we will take ${\bf{y}}(t)=(y_1(t),\dots,y_{n-1}(t))$, ${\bf{x}}(t)=(x_1(t),\dots,x_{n-1}(t))$, and ${\bf{x}}'(t)=(x_1'(t),\dots,x_{n-1}'(t))$ as the moduli of the domain $D_t$ and write ${\bf{M}}(t):=({\bf y}(t),{\bf x}(t),{\bf x}'(t))$. To obtain the conformal equivalence classes from this $3n-3$ dimensional parameter space, we need to identify $({\bf y}(t),{\bf x}(t),{\bf x}'(t))$ and $({\bf{\tilde{y}}}(t),{\bf\tilde{x}}(t),{\bf\tilde{x}}'(t))$, whenever there exists an $a>0$ such that ${\bf y}=a\bf{\tilde{y}}$, $\bf{x}=a\bf{\tilde{x}}$, and ${\bf x'}=a{\bf \tilde{x}'}$. The extra parameter $\bf M$ keeps track of will be reflected in a symmetry (invariance) of the moduli diffusion. For a standard domain the marked points are 0 and $\infty$. For a point $\bf M$ in the ``moduli space" we denote by $D=D({\boldsymbol M})$ the corresponding standard domain.

By boundary correspondence, if $z\in C_j$, then $g_t(z)\in C_j(t)$ and \[
	\Im(g_t(z))=y_j(t).
\]
Thus, by considering the imaginary part of the chordal Loewner equation,
\begin{equation}\label{E:imaginaryCL}
	\partial_t y_j(t)=-\Im(\Psi_t(g_t(z),\xi(t))).
\end{equation}
Further, if 
\[
	z_j(t)=x_j(t)+iy_j(t),\quad z_j'(t)=x_j'(t)+iy_j(t)
\]
are the endpoints of the slit $C_j(t)$, then
\[
	z_j(t)=g_t(\eta_j(t)+iy_j(0)),\quad z_j'(t)=g_t(\eta_j'(t)+iy_j(0)),
\]
where $x_j(0)<\eta_j(t),\eta_j'(t)<x_j'(0)$. Indeed, the pre-images of the tips of $C_j(t)$, that is $\eta_j(t)+iy_j(0)$ and $\eta_j'(t)+iy_j(0)$, are the solutions to the equation
\[
	\frac{\partial}{\partial z}g_t(z)=0,
\]
on the set of prime-ends corresponding to $ C_j\backslash\{z_j(0),z_j'(0)\}$. A tip of $C_j(t)$ cannot be the image of a tip of $C_j$ because then the analytic function $\partial g_t/\partial z$ would not have the required number of zeroes, $2n-2$. 	
	
\begin{lemma}[Motion of moduli---chordal case]\label{L:mm}
The moduli 												\[															{\bf M}(t)=({\bf y}(t),{\bf x}(t),
	{\bf x}'(t))
\] 
satisfy the system of equations
\begin{align}\label{E:mm}
	\partial_t y_j(t)&=\left[\boldsymbol{P}_t^{-1}
	\frac{\partial\boldsymbol{\omega}_t(\xi(t))}{\partial n}\right]_j,
	\notag\\
	\partial_t x_j(t)&=-\Re\left(\Psi_t\left(x_j(t)+iy_j(t),\xi(t)\right)\right),
	\notag\\
	\partial_t x_j'(t)&=
	-\Re\left(\Psi_t\left(x_j'(t)+iy_j(t),\xi(t)\right)\right),
	\end{align}
for  $j=1,\dots,n-1$.
\end{lemma}

\begin{proof}
We note that $\partial g_t/\partial z$ and $\partial^2 g_t/(\partial z)^2$ are analytic functions that extend analytically to the prime-ends corresponding to $C_1,\dots, C_{n-1}$ with the endpoints of the slits removed. By the implicit function theorem, 
\[
	t\mapsto \eta_j(t)+iy_j(0)
\]
is differentiable with derivative 
\[
	DER_t:=\left[
	\frac{\partial^2 g_t}{(\partial z)^2}(\eta_j(t)+iy_j(0))
	\right]^{-1}
	\frac{\partial^2 g_t}{\partial t\partial z}(\eta_j(t)+iy_j(0)).
\]
By counting zeroes we find that 
\[
	\frac{\partial^2 g_t}{(\partial z)^2}(\eta_j(t)+iy_j(t))\neq0
\]
and so $DER_t$ is finite. Hence
\begin{align}
	\partial_t x_j(t)&=\partial_t\Re(g_t(\eta_j(t)+iy_j(0)))\notag\\
	&=-\Re\left(\Psi_t\left(x_j(t)+iy_j(t),\xi(t)\right)\right)
	+\Re\left((\partial_z g_t)(z_j(t)) \times DER_t\right)\notag\\
	&=-\Re\left(\Psi_t\left(x_j(t)+iy_j(t),\xi(t)\right)\right).
\end{align}
In a similar way we obtain the derivative of $x_j'(t)$.
It remains to check that \eqref{E:imaginaryCL} agrees with the first equation in \eqref{E:mm}. To this end we note that
\[
	\Im(\Psi_t(z,\zeta))=\frac{\partial G_t(\zeta,z)}{\partial n_1}
	+\boldsymbol{\omega}_t(z)^T\boldsymbol{P}_t^{-1}
	\frac{\partial\boldsymbol{\omega}_t(\zeta)}{\partial n}.
\]
From the boundary behavior of the Green function and the harmonic measures, it follows that for $z\in C_j(t)$
\[
	\frac{\partial G(\zeta,z;t)}{\partial n_1}=0, 		
	\quad\text{and }\omega_k(z)=\delta_{jk}.
\]
The lemma follows.
\end{proof}

We now have our main existence statement.
	
\begin{theorem}\label{T:Cmm}
Given a continuous function $t\in[0,\infty)\mapsto\xi(t)\in \mathbb R$ 	and the moduli $\bf{M}$ of a standard domain $D$, there exists a unique solution ${\bf M}(t)$ to the 	system \eqref{E:mm} on an interval $[0,t_{\xi})$ with ${\bf M}(0)=\bf{M}$, and where $t_{\xi}$ is 	characterized by 
\[
	t_{\xi}=\inf\{\tau:\lim_{t\nearrow \tau}y_j(t)=0\text{ 
	for some }j\in\{1,\dots, n-1\}\}.
\] 
Further, if $D_t$ is the standard domain determined by ${\bf M}(t)$, and if $\Psi_t(z,\zeta)$ is the holomorphic vector field associated to $D_t$ (cf. Section \ref{ss:Green}), then, for any $z\in D$, the equation
\[	
	\partial_t g_t^{D}(z)=\Psi_t(g_t^D(z),\xi(t)),\quad g_0^D(z)=z,
\]
has a unique solution on $[0,t_z)$, where 
\[
	t_z=\sup\{t\le t_{\xi}:\inf_{s\in[0,t]}|g_s^D(z)-\xi(s)|>0\}.
\]
Finally, for $t<t_{\xi}$ set $K_t=\{z\in D:t_z\le t\}$. Then $g_t^D$ is the canonical conformal map from $D\backslash K_t$ onto $D_t$ with hydrodynamic normalization at infinity.
\end{theorem}

\begin{proof}
For the existence of the solution to the moduli equations \eqref{E:mm} on $[0,t_\xi)$ we need to know that the vector field in \eqref{E:mm} is Lipschitz as a function of ${\bf M}$, with a Lipschitz constant that only depends on distance to $\xi(t)$ of the slit (or slits) nearest to $\xi(t)$. Let $\bf M$ and $\bf{\tilde{ M}}$ be two points in moduli space with corresponding standard domains $D$ and $\tilde{D}$, such that 
\[
	|y_j-\tilde{y}_j|,|x_j-\tilde{x}_j|,| x_j'-\tilde{x}_j'|<\epsilon.
\]
We assume that $\epsilon$ is so small that 
\[
	C_j\cap \tilde{C}_k=\emptyset,\text{ whenever }j\neq k.
\]
Denote $z_j, z_j'$ the endpoints of the slit $C_j$ and $\tilde{z}_j, \tilde{z}_j'$ the corresponding endpoints of $\tilde{C}_j$. Denote $\Psi$ the canonical map for $D$ and $\tilde{\Psi}$ the canonical map for $\tilde{D}$. Then we need to show that 
\begin{equation}\label{E:bigO}
	\tilde{\Psi}(\tilde{z_j})-\Psi(z_j), \tilde{\Psi}(\tilde{z_j}')-\Psi(z_j')
	=O(\epsilon),\quad j=1,\dots,n-1.
\end{equation}
This can be shown as in the radial case by the use of an interior variation that induces a smooth mapping $z\mapsto \tilde{z}$ from $D$ to $\tilde{D}$ which maps slit-endpoints to corresponding slit-endpoints, see \cite{bauer.friedrich:2004}. The non-compactness of the upper half-plane is of no concern as the mapping from $D$ to $\tilde{D}$ may be assumed to be the identity outside of a compact.

The second part of the theorem now follows from general results about ordinary differential equations, exactly as in the simply connected case.
\end{proof}


\subsection{Bilateral case}
As we mentioned before, the bilateral case is similar to the radial case. The geometric description of a bilateral standard domain with $n$ boundary components requires $1+3(n-2)$ real parameters: one for the radius $Q$ of the inner circle, and three for each concentric circular slit. If $C_j$ is one of the interior slits, then $C_j=\{r_j e^{i\theta}, \theta_j\le\theta\le\theta_j'\}$, and we will take $m_j=\ln r_j$, and $\theta_j, \theta_j'$ as parameters to identify $C_j$.

If, in an arbitrary $n$-connected domain $D$, where $n\ge2$, we choose a boundary point $w$ and a boundary component that does not contain $w$, then there is a unique conformal map from $D$ onto a bilateral standard domain, which sends $w$ to 1, and the other distinguished boundary component to the inner boundary circle of the standard domain. Thus the conformal equivalence classes of $n$-connected domains with one marked boundary point and one distinguished boundary component which does not include the marked point are given by  $1+3(n-2)=3n-5$ parameters. We call the parameters 
\[
	(\ln Q, m_2,\dots,m_{n-1},\theta_2,\dots,\theta_{n-1},
	\theta_2',\dots,\theta_{n-1}')
\]
the {\it moduli} of the domain. Note that, unlike in the chordal case, these are true moduli, in the sense that different sets of parameters correspond to different conformal equivalence classes.  

In the bilateral case it was natural to choose the parameter $t=\ln Q$ as time. For a bilateral standard domain $D_t$, where $t=\ln Q$, we let 
\[
	{\bf M}(t)=(m_2(t),\dots, m_{n-1}(t),\theta_2(t),\dots,\theta_{n-1}(t),
	\theta_2'(t),\dots,\theta_{n-1}'(t)).
\]
We then can obtain the following results in the same way as in the chordal case.

\begin{lemma}[Motion of moduli---bilateral case]
The moduli ${\bf M}(t)$ satisfy the system 
\begin{align}
	\partial_t m_j(t)&=1-\Im[\Psi_t(m_j(t)e^{i\theta_j(t)},\xi(t))
	-\Psi_t(e^t,\xi(t))], \notag\\
	\partial_t\theta_j(t)&=\Re[\Psi_t(m_j(t)e^{i\theta_j(t)},\xi(t))
	-\Psi_t(e^t,\xi(t))], \notag\\
	\partial_t\theta_j'(t)&=\Re[\Psi_t(m_j(t)e^{i\theta_j'(t)},\xi(t))
	-\Psi_t(e^t,\xi(t))], 
\end{align}
where $j=2,\dots, n-1$.
\end{lemma}

As in the radial case, it can be shown that the vector field appearing on the right above is Lipschitz in the moduli and we obtain

\begin{theorem}\label{T:Bmm}
Given a continuous function $t\in[0,\infty)\mapsto\xi(t)\in S^1$ 	and the moduli $\bf{M}$ of a bilateral standard domain $D$ with interior boundary circle of radius $Q$, there exists a unique solution ${\bf M}(t)$ to the system \eqref{E:mm} on an interval $[\ln Q,t_{\xi})$ with ${\bf M}(0)=\bf{M}$, and where $t_{\xi}$ is characterized by 
\[
	t_{\xi}=\inf\{\tau:\lim_{t\nearrow \tau}m_j(t)=0\text{ 
	for some }j\in\{2,\dots, n-1\}\}.
\] 
Further, if $D_t$ is the bilateral standard domain determined by ${\bf M}(t)$, and if $\Psi_t(z,\zeta)$ is the holomorphic vector field associated to $D_t$ (cf. Section \ref{ss:Green}), then, for any $z\in D$, the equation
\[	
	\partial_t \ln g_t^{D}(z)=1+[\Psi_t(g_t^D(z),\xi(t))-\Psi_t(e^t,\xi(t))],\quad g_{\ln Q}^D(z)=z,
\]
has a unique solution on $[\ln Q,t_z)$, where 
\[
	t_z=\sup\{t\le t_{\xi}:\inf_{s\in[\ln Q,t]}|g_s^D(z)-\xi(s)|>0\}.
\]
Finally, for $t<t_{\xi}$ set $K_t=\{z\in D:t_z\le t\}$. Then $g_t^D$ is the canonical conformal map from $D\backslash K_t$ onto $D_t$ with $g_t^D(Q)=e^t$.
\end{theorem}


\section{Chordal SLE in multiply connected domains}\label{S:CSLE}

\subsection{Conformal Invariance and Markovian-type Property}
The purpose of this paper is 1) to give a ``natural" construction of conformally invariant measures on ``simple curves" in multiply connected domains, and 2) to study some of the properties of these random curves. We will now motivate, using informal arguments, our particular construction of conformally invariant measures on simple curves. The arguments lead to a small class of processes which contains  chordal $SLE_{\kappa}$ in multiply connected domains.    

For a domain $D$ with $n$ non-degenerate boundary continua and two boundary points (or, more generally, prime ends) $z$ and $w$ lying on the same boundary continuum, let  $W(D,z,w)$ be the set of Jordan arcs in $D$ with endpoints $z$ and $w$. Denote $\{{\mathcal L}_{D,z,w}^{\bf M}\}_{D,z,w}$ a family of probability measures on Jordan arcs in the complex plane such that 
\[
	{\mathcal L}_{D,z,w}^{\bf M}(W(D,z,w))=1,
\]
and where ${\bf M}=M(D)$.
Such families arise, or are conjectured to arise, as distributions of interfaces of statistical mechanical systems at criticality. Based on these models, e.g. percolation, one expects that the distributions describing the interfaces in different domains with different marked points are related by a Markovian-type property and conformal invariance. Denote $\gamma$ a random Jordan arc with law $\mathcal L_{D,z,w}^{\bf M}$. The Markovian-type property says that if $\gamma'$ is a sub-arc of $\gamma$ which has $z$ as one endpoint and whose other endpoint we denote by $z'$, and if ${\bf M}'=M(D\backslash\gamma')$, then the conditional law of $\gamma$ given $\gamma'$ is 
\begin{equation}\label{E:markov}
	\text{law}(\gamma|\gamma')=\mathcal L_{D\backslash\gamma',z',w}^{{\bf M}'}.
\end{equation}
Conformal invariance means that if $f:D\to D'$ is conformal, $z'=f(z)$, $w'=f(w)$, then
\begin{equation}\label{E:ci}
	\mathcal L_{D',z',w'}^{\bf M}
	=f_{*}\mathcal L_{D,z,w}^{\bf M}.
\end{equation}
If \eqref{E:ci} holds, then to understand the family $\{ \mathcal L_{D,z,w}^{\bf M}\}$ it is enough to consider standard domains $D$, take $w=\infty$, $z=0$, and, by the identification of standard domains with their moduli, we may write 
\[
	\mathcal L_{D,0,\infty}^{\bf M}
	=\mathcal L^{\bf M}.
\]
In this case there is a natural parametrization of the Jordan arcs we consider. Let 
\[
	s\in[0,\infty)\mapsto\gamma(s)\in\overline{D}
\]
be a Jordan arc in a standard domain $D$ such that 
\[
	\gamma(0)\in \mathbb R,\ \gamma(0,\infty)\subset D,
	\text{ and }\lim_{t\to\infty}\gamma(t)=\infty.
\]
Denote ${\bf M}=M(D)$ the point in the moduli space corresponding to $D$ and let $g_t^{\bf M}$ be the canonical mapping from $D\backslash\gamma[0,t]$ onto a standard domain $D_t:=g_t^{\bf M}(D\backslash\gamma[0,t])$. Then
\begin{equation}\label{E:expandthis}
	g_t^{\bf M}(z)=z+\frac{a_t}{z}+o(1/|z|),\quad z\to\infty,
\end{equation}
where $a_t$ is called the half-plane capacity. The function $t\mapsto a_t$ is continuous, strictly increasing, starts at zero and satisfies
$a_t\to\infty$ as $t\to\infty$ (this final statement is not true if the curve  creeps along to infinity very close to the real axis and we exclude this case for the purpose of this argument). Thus we may and always will assume that $\gamma$ is parametrized by half-plane capacity, i.e. so that $a_t=2t$. This parametrization is natural in the following sense. If $t\ge0$, ${\bf M}(t)=M(D_t)$, and $\tilde{\gamma}$ is the curve defined by 
\[
	s\in[0,\infty)\mapsto\tilde{\gamma}(s)
	=g_t^{\bf M}(\gamma(t+s)),
\]
then the canonical mapping $g_s^{{\bf M}(t)}$ from $D_t\backslash\tilde{\gamma}[0,s]$ is given by
\[
	g_s^{{\bf M}(t)}=g_{t+s}^{\bf M}\circ
	\left(g_t^{\bf M}\right)^{-1},
\]
and so $g_s^{{\bf M}(t)}(D_t\backslash\tilde{\gamma}[0,s])=D_{t+s}$. In particular, it is easy to see that 
\[
	g_s^{{\bf M}(t)}(z)=z+\frac{2s}{z}+o(1/|z|),\quad z\to\infty, 
\]
i.e. $\tilde{\gamma}$ is also parametrized by half-plane capacity. 

Let now $\{g_s^{\bf M}:s\ge0\}$ be the random family of canonical maps corresponding to the random Jordan arcs $\{\gamma[0,s]:s\ge0\}$ in a standard domain $D$, and denote
\[
	\mathcal{L}^{\bf M}
	=\text{law}(\{g_s^{\bf M}:s\ge0\}).
\]
Then, applying first the Markovian-type property and then conformal invariance, \eqref{E:markov}, \eqref{E:ci}, we find
\[
	\text{law}(\{g_{t+s}^{\bf M}:s\ge0\}|
	g_t^{\bf M})
	=\left(g_t^{\bf M}\right)_*^{-1}
	\mathcal L^{{\bf M}(t)}. 
\]
Equivalently,
\begin{equation}\label{E:law1}
	\text{law}\left(\{g_{t+s}^{\bf M}\circ
	\left(g_t^{\bf M}\right)^{-1}:s\ge0\}|
	g_t^{\bf M}\right)
	=\text{law}(\{g_{s}^{{\bf M}(t)}:s\ge0\}).
\end{equation}
By the chordal Loewner equation, \eqref{E:CL}, for each $t\ge0$, the $\sigma$-field generated by $g_t^{\bf M}$ is equal to $\sigma((\xi(r),{\bf M}(r)):r\in[0,t])$, where $\xi(0)=0$. Similarly, it is easy to see that we can reconstruct $g_{t+s}^{{\bf M}}\circ \left(g_t^{\bf M}\right)^{-1}$ from $\{(\xi(t+r)-\xi(t),{\bf M}(t+r)):r\in[0,s]\}$. Thus \eqref{E:law1} implies
\begin{align}\label{E:m}
	\text{law}&(\{(\xi(t+s)-\xi(t),{\bf M}(t+s)):s\ge0\}|
	\{(\xi(r),{\bf M}(r)):r\in[0,t]\})\notag\\
	&=\text{law}(\{(\tilde{\xi}(s),\tilde{{\bf M}}(s)):s\ge0\}),
\end{align}
where $\tilde{{\bf M}}(s)=M(D_t\backslash\tilde{\gamma}[0,s])$, for a random Jordan arc $\tilde{\gamma}$ with law $\mathcal L^{{\bf M}(t)}$. The equality \eqref{E:m} is precisely the statement that $\{(\xi(t),{\bf M}(t)):t\ge0\}$ is a Markov process. We note that in the simply connected case ($n=1$), \eqref{E:m} reduces to 
\[
	 \text{law}(\{\xi(t+s)-\xi(t):s\ge0\}|
	\{\xi(r):r\in[0,t]\})
	=\text{law}(\{\tilde{\xi}(s):s\ge0\}),
\]
from which it follows that $\xi$ is a process with independent, and identically distributed increments. From this, continuity, and the symmetry $\text{law}(\xi)=\text{law}(-\xi)$, Schramm derived in \cite{schramm:2000} that $\xi(t)=\sqrt{\kappa}B_t$ for a standard one-dimensional Brownian motion and a positive constant $\kappa$. The continuity follows from the continuity of the Jordan arcs, and the symmetry is actually observed in various discrete models, such as the percolation exploration process.


\subsection{Scaling}
For chordal SLE in the upper half-plane $\mathbb H$ the scaling property is usually arrived at as a consequence of the scaling property of the driving function, Brownian motion. Indeed, denote
\[
	\partial_t g_t(z)=\frac{2}{g_t(z)-\sqrt{\kappa}B_t},\quad g_0(z)=z,
\]
chordal SLE in $\mathbb H$ and let $K_t$ be its hull at time $t$, i.e. $g_t$ maps $\mathbb H\backslash K_t$ conformally onto $\mathbb H$. If $c>0$, then $h_t$ defined by
\[
	h_t(z)=\frac{1}{c}g_{c^2t}(cz)
\]
is the normalized conformal map from $\mathbb H\backslash\frac{1}{c} K_{c^2t}$ onto $\mathbb H$ and
\[
	\partial_t h_t(z)=\frac{2}{h_t(z)-\sqrt{\kappa}\frac{1}{c}B_{c^2t}},
	\quad h_0(z)=z.
\]
Since $\frac{1}{c} B_{c^2t}$ is also a standard Brownian motion, it follows that 
\begin{equation}\label{E:scaleconform}
	\text{law}\left(\frac{1}{c}K_{c^2t}:t\ge0\right)
	=\text{law}(K_t:t\ge0).
\end{equation}
However, we can also turn the argument around and ask for a law on growing compacts $K_t$ in the upper half-plane which is conformally invariant, the parameter $t$ being the half-plane capacity as above. For the conformal map $z\mapsto cz$, this implies \eqref{E:scaleconform}, as the half-plane capacity scales quadratically. Denote $\gamma_t$ the tip of the curve generating $K_t$. Then the driving function for the Loewner equation is given by $w_t=g_t(\gamma_t)$, and \eqref{E:scaleconform} implies
\[
	\text{law}\left(\frac{1}{c}w_{c^2t}:t\ge0\right)
	=\text{law}(w_t:t\ge0),
\]
i.e. the driving function has Brownian scaling. Examples of diffusion processes with Brownian scaling are multiples of Brownian motion but also Bessel processes. More generally, if $w$ satisfies the stochastic differential equation
\[
	dw_t=\sigma(w_t)\ dB_t+b(w_t)\ dt,
\]
then $w$ has Brownian scaling if 
\[
	\sigma(cx)=\sigma(x),\quad c b(cx)= b(x),
\]
see \cite{revuz.yor:1999}. If we assume that the coefficients $\sigma$ and $b$ are continuous, then this is saying that $\sigma$ is constant, and $b$ homogeneous of degree minus one. 

In the multiply connected case we can argue similarly. Denote $g_t^{\bf M}$ the normalized conformal map from $D({\bf M})\backslash K_t^{\bf M}$ onto $D_t$. The superscript ${\bf M}$ indicates that the random compact set is a hull in the domain $D({\bf M})$. Conformal invariance of the growing random compacts $K_t^{\bf M}$ requires that 
\begin{equation}\label{E:scaleK}
	\text{law}\left(\frac{1}{c}K_{c^2t}^{c{\bf M}}:t\ge0\right)
	=\text{law}(K_t^{\bf M}:t\ge0).
\end{equation}
Let $w_t^{\bf M}=g_t^{\bf M}(\gamma_t)$, where $\gamma_t$ is the tip of the curve generating $K_t^{\bf M}$. Then \eqref{E:scaleK} implies
\[
	 \text{law}\left(\frac{1}{c}w_{c^2t}^{c{\bf M}}, 
	 \frac{1}{c}{\bf M}_{c^2t}^{c \bf M} :t\ge0\right)
	=\text{law}(w_t^{\bf M}, {\bf M}_t^{\bf M}:t\ge0),
\]
where the superscript ${\bf M}$ indicates that ${\bf M}_0=\bf M$. Thus, the moduli diffusion $(w_t,{\bf M}_t)$ also satisfies Brownian scaling. As in the one dimensional (simply connected) case, this implies under mild regularity assumptions  that the coefficients of the martingale parts of the stochastic differential equation are constant, and the drift coefficients all homogeneous of degree minus one. The drift coefficients of $d{\bf M}_t$ are given in \eqref{E:mm} and we check immediately that they are indeed homogeneous of degree minus one.  


\subsection{Moduli diffusion and interactions with the boundary}
For the purposes of this subsection a different normalization of mappings on standard domains is useful. We will change the normalization of the maps $g_t$ by changing  the vector field in the chordal Loewner equation \eqref{E:CL}. For a chordal standard domain $D$ and $w\in\mathbb R$, define the real function $k(w)$ by
\begin{equation}\label{E:defk}
	k(w)=\lim_{z \to w}(\Psi(z,w)+\frac{2}{z-w}),
\end{equation}
and the conformal map $\Psi^0(z)=\Psi^0(z,w)$ by
\[
	\Psi(z,w)=\Psi^0(z,w)+k(w).
\]
Then $\Psi^0(z,w)=\Psi_D^0(z,w)$ is the unique conformal map from $D$ onto the upper half-plane with a finite number of horizontal slits which sends $w$ to $\infty$ and satisfies
\[
	\lim_{z\to w}(\Psi^0(z,w)+\frac{2}{z-w})=0.
\] 
Consider the modified chordal Loewner equation
\begin{equation}\label{E:mCL}
	\partial_t g_t^0(z)=-\Psi^0(g_t^0(z),\xi^0(t)),\quad g_0^0(z)=z.
\end{equation}
This is the normalization used in \cite{zhan:2004}. Geometrically,  this normalization means that if $g^0$ removes a small vertical slit from the boundary of the upper half-plane, then the images of the two sides of this slit under $g^0$ have the same length up to first order, see \cite{BF:2004a}.   

Let $\kappa$ be a positive real number and $A=A_{\kappa}(w,{\bf M})$  a function homogeneous of degree minus one in the variables $w\in\mathbb R$, and ${\bf M}$ in an open subset of $\mathbb R^{3n-3}$. Consider the system of stochastic differential equations
\begin{align}\label{E:syst}
	d\xi(t)&=\sqrt{\kappa}d B_t+A_{\kappa}(\xi(t),{\bf M}_t)\ dt,\notag\\
	d y_j(t)&=\Im\left(\Psi_t^0\left(x_j(t)+iy_j(t),\xi(t)\right)\right),
	\notag\\
	d x_j(t)&=\Re\left(\Psi_t^0\left(x_j(t)+iy_j(t),\xi(t)\right)\right),
	\notag\\
	d x_j'(t)&=
	\Re\left(\Psi_t^0\left(x_j'(t)+iy_j(t),\xi(t)\right)\right),\quad
	j=1,\dots,n-1,
\end{align}
where ${\bf M}_t=(y_1(t),\dots, y_{n-1}(t), x_1(t),\dots, x_{n-1}(t), x_1'(t),\dots, x_{n-1}'(t))$. If $A$ is Lipschitz, this system has a unique solution. Then we can solve the modified chordal Loewner equation \eqref{E:mCL} for $(\xi(t),{\bf M}_t)$. Denote $K_t$ the random compact such that $g_t^0$ maps the complement of $K_t$ in $D$ conformally onto the standard domain $D_t$.

We can interpret the term $A$ as an interaction of the random growing compact set $K_t$ with the boundary components, and it may be possible to choose $A$ so that the set $K_t$ will avoid these interior boundary components. A similar situation arises for  $\text{SLE}_{\kappa,\rho}$, see \cite{dubedat:2004}. In that case, a random growing compact set in a simply connected domain interacts with a finite number $n$ of boundary points, the interaction strength at point $j$ being given by a real constant $\rho_j$. Then the driving function for the chordal Loewner equation is given by the diffusion
\begin{align}
	d\nu(t)&=\sqrt{\kappa} dB_t
	+\sum_{j=1}^n\frac{\rho_j}{\nu(t)-Z^j(t)}\ dt\notag\\
	dZ^j(t)&=\frac{2}{Z^j(t)-\nu(t)}\ dt,\quad j=1,\dots,n,
\end{align}
a system with drift coefficients homogeneous of degree minus one similar to \eqref{E:syst}. 
  
There are many possible candidates for the homogeneous function $A(w,{\bf M})$. If it is to be a domain functional of the domain $D=D({\bf M})$, then natural candidates arise from derivatives of the Green function. Indeed, if $G(z,w,{\bf M})$ is the Green function for the domain $D=D({\bf M})$ and $c>0$, then 
\[
	G(z,w,{\bf M})=G(cz,cw, c{\bf M})
\]
by conformal invariance and so 
\[
	\frac{\partial_z^k\partial_w^l 
	G(z,w,{\bf M})}{\partial_z^m\partial_w^n G(z,w,{\bf M})}
\]
is homogeneous of degree minus one whenever
\[
	k+l=m+n+1,\quad k,l,m,n\in\mathbb N.
\]
The ``harmonic random Loewner chains''  Zhan studies in his thesis, see \cite{zhan:2004}, correspond to the choice $k=l=m=1$, $n=0$. Via integration, or directly by conformal invariance, we also see that 
\[
	\frac{\partial_z^{k+1}\omega_j(z,{\bf M})}{\partial_z^k\omega_j(z,{\bf M})}
\]
is homogeneous of degree minus one.  

	
\subsection{Chordal SLE, percolation, and locality}

The case of percolation is an example where there is no interaction, that is $A\equiv0$. For the following calculation we return to the original chordal Loewner equation \eqref{E:CL}. Then $\xi$ in \eqref{E:syst} has a nonzero drift coming from changing back the normalization.Thus, to model cluster-boundaries of percolation in a multiply connected domain $D$ we make the ansatz
\begin{equation}\label{E:ansatz}
	d\xi(t)=-k_t(\xi(t))+\sqrt{\kappa}\ dB_t,
\end{equation}
where the subscript $t$ refers to the domain $D_t$, $k_t$ to \eqref{E:defk}, and where ${\bf M}(t)$ satisfies \eqref{E:mm}.

This choice of drift reflects that the exploration process for percolation is as likely to turn right as it is to turn left. Other discrete models lead to different drifts. 
In this section we show that the ansatz \eqref{E:ansatz} leads to random growing compacts satisfying the locality property if $\kappa=6$. 

Denote $\{g_t^E,t\ge0\}$ the solution of the chordal Loewner equation in a standard domain $E$  starting at $z=0$ for the diffusion \eqref{E:ansatz}. Denote  $\{K_t,t\ge0\}$ the associated growing compacts. Let $A$ be a hull in $E$ that does not contain zero. For the following calculations we restrict to the event $\{t<\tau\}$, where $\tau:=\inf\{t:K_t\cap A\neq\emptyset\}$. Let $\Phi_A$ be the canonical mapping from $E\backslash A$, $g_t^*$ the canonical mapping from $\Phi_A(E\backslash(A\cup K_t))$, and $h_t$ the canonical mapping from  $g_t(E\backslash(A\cup K_t))$. Since the canonical mapping for $E\backslash(A\cup K_t)$ is unique, we have
\begin{equation}
	h_t\circ g_t=g_t^*\circ\Phi_A.
\end{equation}
Furthermore, up to a time change, the family $\{g_t^*\}$ also satisfies a chordal Loewner equation beginning with the standard domain $E^*:=\Phi_A(E\backslash A)$. In fact, reasoning as in \cite{BF:2004a}, it follows that 
\begin{equation} \label{E:KLstar}
	\partial_t g_t^*(z)
	=-|h_t'(\xi(t))|^2\Psi_t^*(\xi^*(t),w_t^*),
\end{equation}
where $w_t^*=g_t^*(z)$, and $\xi^*(t)=h_t(\xi(t))$. The question we are interested in is whether $(\xi^*,{\bf M}^*)$ is a time change of $(\xi,{\bf M})$. Since $h_t=g_t^*\circ\Phi_A\circ g_t^{-1}$, we have
\begin{equation}\label{E:ht}
	\partial_t h_t(z)=\left[\partial_t g_t^*\right](\Phi_A(g_t^{-1}(z)))+(g_t^*\circ\Phi_A)'(g_t^{-1}(z))(\partial_t g_t^{-1}(z)),
\end{equation}
and we note that 
\begin{equation}\label{E:g inverse}
	\partial_t g_t^{-1}(z)=(g_t^{-1})'(z)\Psi_t(\xi(t),z).
\end{equation}
Then \eqref{E:ht},\eqref{E:KLstar}, and \eqref{E:g inverse} imply
\begin{equation}\label{E:phit}
	\partial_t h_t(z)
	=- h_t'(\xi(t))^2\Psi_t^*(\xi^*(t),h_t(z))+h_t'(z)\Psi_t(\xi(t),z).
\end{equation}
Hence the stochastic differential 
\[
	\partial_t h_t(z)\ dt+h_t'(\xi(t))\ d\xi(t)
\]
has martingale part $h_t'(\xi(t))\sqrt{\kappa}\ dB_t$ and its drift part can be grouped into the three components
\begin{align}
	I:&=- h_t'(\xi(t))^2[\Psi_t^*(\xi^*(t),h_t(z))-k_t^*(\xi^*(t))]\ dt\notag\\
	&\quad+h_t'(z)[\Psi_t(\xi(t),z)-k_t(\xi(t))]\ dt,\notag\\
	II:&=-h_t'(\xi(t))^2 k_t^*(\xi^*(t))\ dt,\notag\\
	III:&=-[h_t'(\xi(t))-h_t'(z)]k_t(\xi(t))\ dt.
\end{align}
When $z\to\xi(t)$, then part $III$ converges to zero, and part $II$, together with the martingale part, converges to a time-change  of \eqref{E:ansatz} starting at $E^*$. Finally, for part $I$, by the definition of $k(\xi;t)$ a double application of l'H\^opital's rule gives 
\begin{equation}
	\lim_{z\to \xi}\left(\frac{2 h'(\xi)^2}{h(z)-h(\xi)}-
	\frac{2 h'(z)}{z-\xi}\right)=-3 h''(\xi).
\end{equation}
Thus, by It\^o's formula,
\begin{equation}\label{E:ito}
	d h_t(\xi(t))
	=-h_t'(\xi(t))^2 k_t^*(\xi^*(t))\ dt
	+\frac{\kappa-6}{2}h_t''(\xi(t))\ dt+h_t'(\xi(t))\sqrt{\kappa}\ dB_t,
\end{equation}
which is indeed a time-change of \eqref{E:ansatz} if and only if $\kappa=6$. From \eqref{E:KLstar} it follows immediately that the equations for ${\bf M}^*$ are given by the same time change of the equations for ${\bf M}$.  

\begin{theorem}[Chordal $\text{SLE}_6$]
The solution to the chordal Loewner equation based on the diffusion \eqref{E:ansatz} satisfies the locality property if and only if $\kappa=6$.
\end{theorem}



\begin{thebibliography}{99}

\bibitem{ahlfors:1966}
 	L. Ahlfors, {\it Complex Analysis}, 2nd. ed., McGraw-Hill, New 		York, 1966.

\bibitem {Ai}
M. Aizenman,
{\it The geometry of critical percolation and conformal invariance}, Stat. Phys.  {\bf 19} (1996), 104--120.


\bibitem{bauer.friedrich:2004}
	R. Bauer, R. Friedrich, {\it Stochastic Loewner evolution in multiply connected domain},
C. R. Acad. Sci. Paris, Ser. I {\bf 339}, 579-584 (2004).

\bibitem{BF:2004a}
	R. Bauer, R. Friedrich, {\it On radial stochastic Loewner evolution in multiply connected domains},
	arXiv. 


\bibitem{BPZ}
{A.A. Belavin, A.M. Polyakov, A.B. Zamolodchikov, {\it 
Infinite conformal symmetry in two-dimensional quantum field theory},
Nuclear Phys. B {\bf 241} (1984), 333--380.}



\bibitem{cardy:2004}
J. Cardy, {\it SLE(kappa,rho) and Conformal Field Theory}, arXiv
math-ph/0412033.

\bibitem{courant:1950}
	R. Courant, {\it Dirichlet's Principle}, with an appendix by M. 	Schiffer, Interscience, New York, 1950.

\bibitem{dubedat:2004}
	J. Dubedat, {\it Some remarks on commutation relations for 		SLE}, arXiv, math.PR/0411299.

\bibitem {FK}
{R. Friedrich, J. Kalkkinen, 
On conformal field theory and stochastic Loewner evolution,
Pr\'epublications de l'IH\'ES~P/03/28, (2003), which appeared inÊ 
Nuclear Phys. B {\bf 687} (2004), no. 3, 279--302.}
	
\bibitem{komatu:1943}
  Y. Komatu, {\it Untersuchungen \"uber konforme {A}bbildung von zweifach zusammenh\"angenden {G}ebieten},   Proc. Phys.-Math. Soc. Japan (3) {\bf 25} (1943), 1--42.

\bibitem{komatu:1950}
    Y. Komatu, {\it On conformal slit mapping of multiply-connected domains}, Proc. Japan Acad. {\bf 26} (1950), no.~7, 26--31.
    
\bibitem{KBonn} 
{M.~Kontsevich,
Arbeitstagung 2003 - ``CFT, SLE and phase boundaries"
MPIM (2003).}

\bibitem {LPS}
R. Langlands, Y. Pouillot, Y. Saint-Aubin, {\it
 Conformal invariance
        in two-dimensional percolation}, Bull. A.M.S. {\bf 30} (1994), 1--61.


\bibitem{lsw:2001a}
	G. F. Lawler, O. Schramm\ and\ W. Werner, {\it Values of Brownian intersection exponents. I. Half-plane exponents}, Acta Math. {\bf 187} (2001), no.~2, 275--308.


\bibitem{lsw:2003}
G. F. Lawler, O. Schramm\ and\ W. Werner, {\it Conformal restriction: the chordal case}, J. Amer. Math. Soc. {\bf 16} (2003), 917--955.

\bibitem{lsw:2004}
	G. F. Lawler, O. Schramm\ and\ W. Werner, {\it Conformal invariance of planar loop-erased random walks and uniform spanning trees}, Ann. Probab. {\bf 32} (2004), no.~1B, 939--995.

\bibitem{loewner:1923}
	K. L\"owner, {\it Untersuchungen \"uber schlichte konforme Abbildungen des Einheitskreises I.}, Math. Ann. {\bf 89} (1923), 103--121.


\bibitem{nehari:1952}
	Z. Nehari, {\it Conformal Mapping}, McGraw-Hill, New York, 1952.

\bibitem{pommerenke:1992}
	C.ÊPommerenke, {\it Boundary behaviour of conformal maps}, Grundlehren der Mathematischen Wissenschaften, 299, Springer-Verlag, Berlin, 1992.

\bibitem{revuz.yor:1999}
	D. Revuz, M. Yor, {\it Continuous Martingales and Brownian Motion}, Grundlehren der mathematischen Wissenschaften, Vol. 293, 3rd edition, Springer, Heidelberg, 1999.

\bibitem{rohde.schramm:2003}
	S. Rohde, O. Schramm, {\it Basic properties of SLE}, preprint, 		arXiv:math.PR/0106036 v2.

\bibitem{schiffer:1946}
  M. Schiffer, {\it Hadamard's formula and variation of domain-functions}
  Amer. J. Math. {\bf 68} (1946), 417--448.
  
 \bibitem{schiffer.spencer:1954}
 	M. Schiffer, D. Spencer, {\it Functionals of finite Riemann surfaces}, Princeton University Press, Princeton, New Jersey, 1954.

\bibitem{schramm:2000}
	O. Schramm, {\it Scaling limits of loop-erased random walks 	and uniform spanning trees}, Israel J. Math. {\bf 118} (2000), 	221--288.


\bibitem{werner:2003}
	W. Werner, {\it Random planar curves and Schramm-Loewner evolutions}, lecture notes from the 2002 St. Flour summer school, Springer, Berlin, 2003.

\bibitem{zhan:2004}
	D. Zhan, {\it Random Loewner chains in Riemann surfaces}, thesis, California Institute of Technology, 2004.

\end{thebibliography}
\end{document}